\documentclass[12 pt, a4paper]{article}
\usepackage{graphicx} 
\usepackage{amsfonts,amsmath,amssymb,amsthm,url,float,enumitem}
\usepackage[hmargin=1.1in,vmargin=1.1in]{geometry}
\usepackage{enumitem}

\newcommand{\ignore}[1]{}

\newtheorem{theorem}{Theorem}[section]
\newtheorem{claim}[theorem]{Claim}

\newtheorem{lemma}[theorem]{Lemma}
\newtheorem{corollary}[theorem]{Corollary}

\newcommand{\pts}{\mathcal P}
\newcommand{\curves}{\Gamma}
\newcommand{\RR}{\ensuremath{\mathbb R}}
\newcommand{\R}{\ensuremath{\mathbb R}}
\newcommand{\Z}{\ensuremath{\mathbb Z}}

\def\B{\mathcal{B}}
\def\M{\mathcal{M}}

\def\eps{\varepsilon}

\newcommand{\parag}[1]{\vspace{2mm}

\noindent{\bf #1} }

\usepackage{xcolor}
\usepackage{tikz}
\usetikzlibrary{positioning, 
                quotes}

\title{Distinct Distances under well-behaved norms\thanks{This work was done as part of the 2026 NYC Discrete Math REU, supported by NSF award DMS-2349366 and by Jane Street.}}
\author{Lillian MacArthur\thanks{Harvard University, {\sl lcmacarthur@college.harvard.edu}}
\and
Adam Sheffer\thanks{CUNY, {\sl adamsh@gmail.com}. Supported by PSC-CUNY award GR-00017710.}}
\date{}

\begin{document}

\maketitle

\begin{abstract} 

We improve our understanding of the Erd\H{o}s' distinct distances problem under non-Euclidean norms, in two ways: \\[2mm]
(i) For norms that are definable in an $o$-minimal structure and whose non-line bisectors have a finite intersection, we obtain a lower bound of $\Omega(n^{5/6})$ for the minimum number of distinct distances. 
In particular, this improves the current best bound of $\Omega(n^{4/5})$ for all $\ell_p$ norms with non-integer $p>1$. \\[2mm]
(ii) We derive the stronger bound $\Omega(n^{6/7-\eps})$ under all semi-algebraic smooth norms.

Our proof technique is based on deriving properties of bisectors under various norms.
\end{abstract}

\section{Introduction}

Several central problems of discrete geometry have a misleading quality: The problem can be easily described to a high-school student, creating a first impression of being a riddle. 
However, studying the problem requires developing a deep theory and understanding properties of the underlying geometry. 
One classic example is the \emph{Erd\H{o}s distinct distances problem}.
One aspect of this problem asks what distinguished the Euclidean distance from most other distance norms.  

In 1946, Erd\H{o}s \cite{Erdos46} asked for the minimum number of distinct distances spanned by any set of $n$ points in $\R^2$. 
Denoting this quantity as $d(n)$, he provided a lower bound of $d(n) = \Omega(n^{1/2})$ and an upper bound of $d(n) = O(\frac{n}{\sqrt{\log n}})$. 
Erd\H{o}s conjectured that $d(n) = \Theta(\frac{n}{\sqrt{\log n}})$.
After about 65 years and a huge body of work dedicated to the problem,
Guth and Katz \cite{GK15} proved that $d(n) = \Omega(\frac{n}{\log n})$.
While a gap of $\sqrt{\log n}$ remains, the problem is now considered mostly solved.

Throughout this work, we work in $\RR^2$ without explicitly mentioning this. 
In particular, when discussing norms, these are defined in $\RR^2$.
For a technical discussion of norms and the distances associated with them, see Section \ref{sec:prelim}.
Alon, Buci\'c, and Sauermann \cite{ABS25} proved the following result, building upon a previous work of Matou\v sek \cite{Matousek11}.

\begin{theorem} \label{th:MostNorms}
For most norms, the minimum number of distances spanned by $n$ points is $n-o(n)$.
\end{theorem}

In Theorem \ref{th:MostNorms}, the notation $o(n)$ means asymptotically smaller than $n$. 
The term ``most norms'' can be thought of topologically: When considering the space of all norms, the subspace of norms that do not satisfy the above is of a lower dimension. 
For more rigorous details, see \cite{ABS25}.

The Euclidean norm is not standard in the sense of Theorem \ref{th:MostNorms}, since Erd\H os's construction has $\Theta(\frac{n}{\sqrt{\log n}})$ distances.
This leads to difficult questions: \emph{What makes the Euclidean norm different from most other norms? 
Are other common norms, such as $\ell_p$ norms, special in the same way?}

Embarrassingly, almost nothing is known about the above questions. 
Not a single specific norm is known to be typical in the sense of Theorem \ref{th:MostNorms}.
In addition, not a single non-Euclidean norm is known to be non-typical.\footnote{This statement excludes two uninteresting cases:  ellipsoidal norms and non-strictly convex norms. For more details about these cases, see Section \ref{sec:prelim}.}

Garibaldi studied the distinct distances problem under different norms in her dissertation \cite{Garibaldi04}, supervised by John Garnett and Terence Tao.
She defined a family of well-behaved norms named \emph{potato metrics}, which includes all $\ell_p$ distances.
A potato metric is a strictly convex norm where (i) The unit
circle has finitely many axes of symmetry; (ii) Every two non-line bisectors intersect at a bounded number of
points.\footnote{For definitions of objects such as unit circles and bisectors, see Section \ref{sec:prelim}. Garibaldi first defines potato metrics somewhat differently and then extends the definition to the one stated here.}
The Euclidean norm is the only one that violates condition (i), so we may replace this condition with the statement that the norm is non-Euclidean.
For potato metrics, Garibaldi proved that every $n$ points span $\Omega(n^{4/5})$ distinct distances. 
For $\ell_p$ distances with integer $p$, this bound was improved in \cite{Polymath20} to $\Omega(n^{6/7-\eps})$, for any $\eps>0$.

The current work derives two new distinct distances bounds for non-Euclidean norms.
The first is for norms that are definable in an $o$-minimal structure.
Intuitively, \emph{o-minimality} is a significant generalization of semi-algebraic sets, which still yields tame geometric behavior. 
For example, curves in an $o$-minimal structure intersect in a finite union of points and arcs.
A norm from an $o$-minimal structure may involve exponential functions, logarithmic function, trigonometric functions (under certain restrictions), inverse trigonometric functions, integration such as $y=\int_0^x \frac{e^t}{t}dt$, and much more. 

Examples of curves that are not definable in an $o$-minimal structure: An unbounded sine wave and certain unions of infinitely many arcs.
For additional details about $o$-minimal structures, see Section \ref{sec:prelim}.

For a norm $\|\cdot\|$, let $d_{\|\cdot\|}(n)$ be the minimum number of distances that $n$ points may span under $\|\cdot\|$.
We say that a norm is $o$\emph{-minimal} if its distance equation is definable in an $o$-minimal structure. 

\begin{theorem} \label{th:o-minimal}
    Let $\|\cdot\|$ be a strictly convex $o$-minimal norm for which every two non-line bisectors have a finite intersection. Then, $d_{\|\cdot\|}(n) = \Omega(n^{5/6}).$
\end{theorem}

We may think of the norms of Theorem \ref{th:o-minimal} as the well-behaved potato norms. 
That is, potato norms where unit circles and perpendicular bisectors have nice properties such as bounded intersections. 
One might refer to these as \emph{sweet potato metrics}.

Since $\ell_p$ norms satisfy the assumptions of Theorem \ref{th:o-minimal}, we obtain a new bound for $\ell_p$ norms with non-integer $p$.
The previous best bound, by Garibaldi, is $\Omega(n^{4/5})$.

\begin{corollary} \label{co:lp}
   For all $p > 1,$ we have that $d_{\ell_p}(n) = \Omega(n^{5/6})$.
\end{corollary}

The following result is another useful application of Theorem \ref{th:o-minimal}. 
We say that a curve is \emph{analytic} if it can be described as $y=f(x)$ with analytic $f$.

\begin{corollary} \label{co:o-minimal+analytic}
    Let $\|\cdot\|$ be a strictly convex $o$-minimal norm with analytic bisectors. Then $d_{\|\cdot\|}(n) = \Omega(n^{5/6})$.
\end{corollary}

We say that a norm is \emph{semi-algebraic} if its unit circle is semi-algebraic. 
(For a definition of semi-algebraic sets, see Section \ref{sec:prelim}.) 
Our second distinct distances bound is for semi-algebraic norms.

\begin{theorem} \label{th:semi-alg}
    Let $\|\cdot\|$ be a strictly convex semi-algebraic norm where bisectors intersect finitely many times. Then all $\eps > 0$ satisfy that $d_{\|\cdot\|}(n) = \Omega(n^{6/7 - \epsilon}).$
\end{theorem}

As in the $o$-minimal case, we may apply Theorem \ref{th:semi-alg} for any strictly convex semi-algebraic norm with analytic bisectors. 
The proof of this claim is identical to the proof of Corollary \ref{co:o-minimal+analytic}, so we do not include it here.

Our proof technique is based on studying properties of bisectors under various types of norms. 
For example, we prove the following result.

\begin{theorem} \label{th:LineBisectors}
Let $\|\cdot\|$ be a strictly convex $o$-minimal norm. 
If two pairs of points $(p_1,q_1) \neq (p_2,q_2)$ satisfy $\B(p_1, q_1) = \B(p_2, q_2)$, then this bisector is a line.
\end{theorem}

Despite the above results, we are still far from understanding the behavior of the distinct distances problem under other norms. 
A few examples of central problems that remain unsolved:
\begin{itemize}[noitemsep,topsep=1mm]
\item Find a strictly convex norm $\|\cdot\|$ for which $d_{\|\cdot\|}(n)=\Theta(n)$. (Followup: Find a strictly convex norm $\|\cdot\|$ for which $d_{\|\cdot\|}(n)=n-o(n)$.)
\item Is there a non-ellipsoidal strictly convex norm $\|\cdot\|$ for which $d_{\|\cdot\|}(n)=o(n)$.
\item Find the value of $d_{\ell_p}(n)$ for any finite $p>1$ (or for any other sweet potato norm).
\item Is being definable in an $o$-minimal structure the natural notion for well-behaved norms? If not, should we replace it with a more restrictive definition or with a less restrictive one?
\end{itemize} 

In Section \ref{sec:prelim}, we go over the technical preliminaries in more detail.
For example, we discuss distance metrics and $o$-minimal structures in more detail.
In Section \ref{sec:bisectors}, we study more advanced properties of bisectors. 
In Section \ref{sec:o-min}, we prove our $o$-minimal distance results: Theorem \ref{th:o-minimal} and its applications Corollary \ref{co:lp} and Corollary \ref{co:o-minimal+analytic}. 
In Section \ref{sec:semi-alg}, we prove our semi-algebraic result, Theorem \ref{th:semi-alg}.

\parag{Acknowledgements.}
We thank Aaron Anderson, Nathaniel Nadler, Abhiram Natarajan, and Konrad Swanepoel for helpful conversations. 

\parag{Statement of AI usage.} The authors used ChatGPT (OpenAI) to assist with literature searches, draw some of the figures, and verify the accuracy of some technical statements. 
However, AI was not used to come up with any parts of the proof.  

\section{Technical preliminaries} \label{sec:prelim}

In this section, we discuss technical definitions and tools that are used in our proofs. 
A reader may wish to skim this section and return to it later as needed.

\parag{Norms, metrics, and unit circles.}
Recall that we only consider point sets and norms in $\RR^2$.
A norm $\|\cdot\|$ induces a \emph{metric} or \emph{distance}: We write the distance between points $a,b\in \RR^2$ as $\|p-q\|$.

When it comes to distance problems, it is customary to only consider metrics that are induced by norms. 
While there are metrics that are not induced by norms, these tend not to make sense in our context.
For example, under such metrics, translations may change distances. 
We follow the standard approach, considering only metrics induced by norms. 
We also try to avoid the term ``metrics,'' except when it comes to ``potato metrics.''
 
One might say that the most well-known family of norms are the $\ell_p$ norms. 
For $p\in [1,\infty)$, the $\ell_p$ norm of the vector $(a_x,a_y)$ is $(|a_x|^p+|a_y|^p)^{1/p}$.
This induces the distance
\[ \|a-b\| = \left(|a_x-b_x|^p+|a_y-b_y|^p\right)^{1/p}. \] 
When $p=2$, this is the Euclidean distance.

\begin{figure}[htp]     \centering
    \includegraphics[width=11cm]{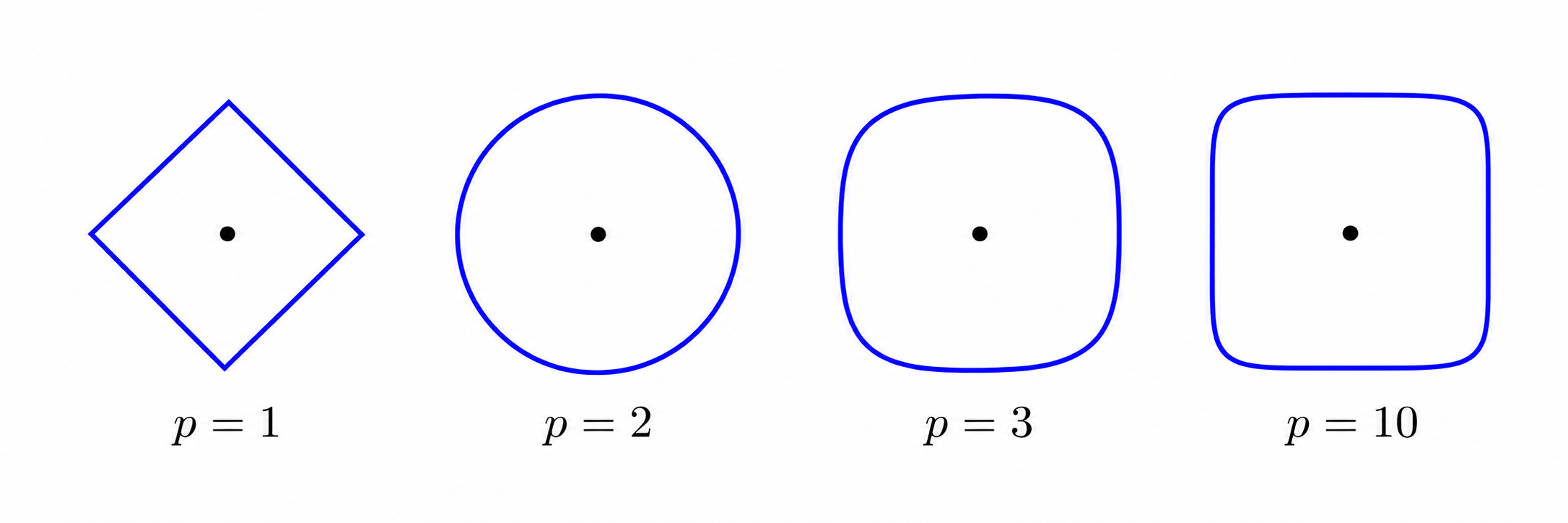}
\vspace{-5mm}    
    
    \caption{Unit circles under several $\ell_p$ norms.}\label{fi:unit}
\end{figure}

The \emph{unit circle} of a norm $\|\cdot\|$ is 
$\{v \in \R^2 : \|v\| = 1\}$.
See Figure \ref{fi:unit}.
By the symmetry of $\|\cdot\|$, the unit circle has 180 degree rotational symmetry around the origin.
Since norms satisfy a triangle inequality, the unit circle is convex. 
We say that a norm $\|\cdot\|$ is \emph{strictly convex} when its unit circle is strictly convex. 
In other words, when the boundary of the unit circle contains no line segments.
This is equivalent to $\|a + b\| < \|a\| + \|b\|$ holding for every linearly independent $a,b$.
While the $\ell_1$ norm is not strictly convex, for every $p\in(1,\infty)$, the $\ell_p$ norm is strictly convex.
See Figure \ref{fi:unit}.

\parag{Special cases: Ellipsoidal norms and non-strictly convex norms.}
We now discuss the two families of norms for which somewhat accurate distinct distances information is known.
An \emph{ellipsoidal norm} is a norm whose unit circle is an ellipse. 

\begin{claim} \label{cl:ellipsoidal}
Let $\|\cdot\|$ an ellipsoidal norm. Then the distinct distances problem with  $\|\cdot\|$ is equivalent to the Euclidean distinct distances problem.
\end{claim}
\begin{proof}[Proof sketch.]
Let $T$ be an invertible linear transformation that takes the circle defined by $y^2+x^2=1$ to the unit circle of $\|\cdot\|$.
Since this unit circle is an ellipse, such a transformation exists.

We note that the Euclidean distance between points $p,q\in\RR^2$ is equal to $\|T(p)-T(q)\|$.
Thus, $T$ is a bijection between point set that maintains the same set of distances (according the Euclidean distance and $\|\cdot\|$, respectively). 
\end{proof}

In a non-strictly convex norm, the unit circle contains a straight-line segment. 
For such a norm, it is not difficult to show that the minimum number of distinct distances is $\Theta(\sqrt{n})$. 
For example, see \cite[Section 2.2]{Garibaldi04}.

\parag{Semi-algebraic sets.}
A set $S\subset \RR^d$ is \emph{semi-algebraic} if it can be expressed as a boolean combination of finitely many polynomial equations and inequalities. 
For example, the half-circle 
\[ \left\{(x,y)\in \RR^2\ :\ x^2+y^2 = 1 \text{ and } x\ge 0\ \right\} \subset \RR^2. \]
is semi-algebraic.
Indeed, the definition is the concatenation of one polynomial equation and one polynomial inequality.
For a deeper discussion of semi-algebraic sets, see for example \cite{BPR06,BCR13,Sheffer2022}.

\parag{$O$-minimal structures.}
For simplicity, we only discuss $o$-minimal structures over $\RR$. 
For additional information, see for example \cite{Coste99,vdDries98}.

In general, an $o$-minimal structure $S$ is a collection of subsets of $\RR^d$ that satisfy certain properties.
Here, $d$ is not fixed, in the sense that elements of $S$ may be subsets of spaces of different dimensions. 
Identical sets in different spaces are considered as distinct elements of $S$.
For example, $\emptyset\subset \RR^2$ and $\emptyset\subset \RR^3$ are distinct elements of $S$.

The elements of $S$ are called \emph{definable sets}.
For $S$ to be an $o$-minimal structure, the following need to hold:
\begin{itemize}[noitemsep,topsep=1mm]
\item For every $d$, all semi-algebraic sets of $\RR^d$ are definable sets, including $\RR^d$ and $\emptyset$. 
\item $S$ is closed under finite intersections, unions, complements, cartesian products, and coordinatewise projections. While intersections, unions, and inverses maintain the dimension of the containing space $\RR^d$, this is not the case for cartesian products and projections. 
For example, a projection from $\RR^d$ that removes a single coordinate leads to an element of $\RR^{d-1}$.
\item The definable sets that are subsets of $\RR^1$ are exactly the finite unions of points and open intervals, including unbounded intervals.
\end{itemize}

The $o$-minimal structure $\RR_{exp}=(\RR; +,\cdot,1,0,exp)$ is the collection of subsets of $\RR^d$ that can be described by first order formulas\footnote{We recall that, \emph{first-order} means that we can only quantify over our real variables. For example, we cannot quantify over functions in these variables. } with real variables, real constants, addition, multiplication, inequalities, equalities, the exponential function $e^x$, the logical operations \emph{and}, \emph{or}, and \emph{not}, and the quantifiers \emph{there exists} and {for every}.
A couple of examples of definable sets of $\RR_{exp}$:
\begin{align*}
&\hspace{6mm}\left\{(x,y)\in\RR^2\ :\ y=x^2\cdot e^{3x^2} \right\}, \\[2mm]
&\left\{x\in R\ :\ \exists y\ (y^2=x \text{ and } e^y<30)\ \right\}.
\end{align*}

A significantly more general $o$-minimal structure is $\RR_\text{Pf}$, the expansion of the real field generated by all \emph{Pfaffian functions}.
The definition of this structure is significantly more involved, so we do not include it here (see for example \cite{Speissegger99}). 
This $o$-minimal structure includes sets defined by combinations of exponential functions, logarithmic functions, trigonometric functions, inverse trigonometric functions, functions with integrals such as $y=\int_0^x \frac{e^t}{t}dt$, and much more. 
(It only supports sines and cosines with bounded ranges.)

From the above examples, one might wonder if $o$-minimal structures cover all sets imaginable. 
This is not the case.
For example, an infinite sine wave is not definable in any $o$-minimal structure.
Neither do various curves that are obtained as the union of infinitely many arcs.

Intuitively, $o$-minimal structures are collections of sets that are well-behaved in several ways.
For example, each set must have finitely many connected components.
Every definable set admits a finite cell decomposition: It can be partitioned into finitely many simple definable pieces, each homeomorphic to an open box. 

A \emph{family of definable sets in} $\RR^2$ is a definable set with $d+2$ variables.
When fixing values for the last $d$ variables, we receive the definiion of a definable set in $\RR^2$. 
For example, the family of circles in $\RR^2$ can be defined with
\[ f(x,y,a,b,r) = (x-a)^2+(y-b)^2-r^2 \]
Indeed, fixing the values of $a,b,r$ leads to a specific circle. 
The \emph{dimension} of a family is the number of parameters $d$ that define a set in $\RR^2$.
For example, the family of circles in $\RR^2$ is three-dimensional.
When the defining formula is semi-algebraic, we may also use the term \emph{semi-algebraic family}.

When studying sets in $\RR^2$ that are defined by polynomials, we have bounds that depend on the degrees of these polynomials. For example, such bounds exist for the number of connected components of an algebraic curve and for the number of intersection points of two algebraic curves with no common components.

A similar situation occurs for sets definable in an $o$-minimal structure. In this setting, the bounds depend on the defining formulas, but not on the values of their parameters. 
That is, the bounds are uniform over all curves from the same definable family. This can be viewed as an analogue of bounded description complexity in the algebraic setting. 
For more information, see for example \cite[Theorem 3.12]{Coste99}.

\begin{lemma} \label{le:components}
Let $F$ be a definable family of subsets of $\RR^2$ in an $o$-minimal structure. Then there exists a constant $N$, depending only on $F$, such that every set in $F$ has at most $N$ connected components.
\end{lemma}

\parag{Geometric incidences.}
Let $\M$ be an $o$-minimal structure.
Let $\pts$ be a set of points and let $\curves$ be a set of definable sets, both in $\RR^2$.
An \emph{incidence} is a pair $(p,\gamma)\in\pts\times\curves$, such that $p\in \gamma$.
Let $I(\pts,\curves)$ denote the number of incidences in $\pts\times\curves$.
The study of incidences is an active research front, with many applications.
For more information about this topic, see for example \cite{Sheffer2022}.

\begin{figure}[htp]     \centering
    \includegraphics[width=10cm]{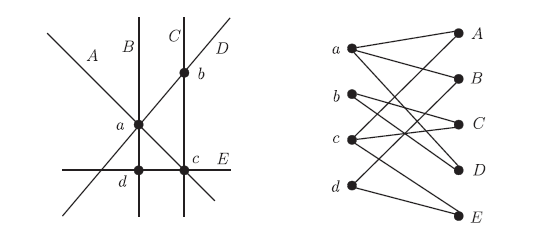}
\vspace{-5mm}    
    
    \caption{A point-line configuration (left) and its incidence graph (right).}\label{fi:IncGraph}
\end{figure}

The \emph{incidence graph} of $\pts\times\curves$ is a bipartite graph, where one side has a vertex for each point of $\pts$ and the other has a vertex for each set of $\curves$. 
There is an edge between a vertex of a point $p$ and a vertex of a set $\gamma$ if $p\in\gamma$.
For an example, see Figure \ref{fi:IncGraph}.
Recall that $K_{s,t}$ is a complete bipartite graph with $s$ vertices on one side and $t$ on the other.
In Figure \ref{fi:IncGraph}, the graph does not contain a $K_{2,2}$.
Indeed, at most one line can pass through any pair of points.

The following is a special case of a result by Chernikov, Galvin, and Starchenko \cite[Theorem 5.14]{CGS20}.

\begin{lemma} \label{ominlemma}
Let $\pts$ be a set of $m$ points in $\RR^2$.
For an $o$-minimal structure $\M$, let $F$ be a $d$-dimensional family of definable sets in $\RR^2$. 
Let $\curves$ be a collection of $n$ sets from $F$, such that the incidence graph of $\pts\times\curves$ contains no $K_{r,r}$.
Then 
\[
    I(\pts,\curves) = O_r \left(m^{\frac{d}{2d - 1}} n^{\frac{2d - 2}{2d - 1}} + m + n \right).
    \]
\end{lemma}

When the sets are semi-algebraic, Sharir and Zahl \cite{SZ17} obtain a stronger incidence bound. 

\begin{lemma} \label{semialglemma}
Let $\pts$ be a set of $m$ points and let $\curves$ be a set of $n$ semi-algebraic sets, both in $\RR^2$.
Assume that the sets of $\curves$ belong to an $s$-dimensional family of constant description complexity and that no two sets of $\curves$ share a common infinite component.  
Then for each $\eps > 0$, we have that 
\[ I(\pts, \curves) = O\left(m^{\frac{2s}{5s - 4}} n^{\frac{5s - 6}{5s -4} + \epsilon} + m^{2/3}n^{2/3} + m + n\right) \]
\end{lemma}

One interesting open problem is to generalize Lemma \ref{semialglemma} to curves defined by Pfaffian functions or, more ambitiously, to curves that are definable in an $o$-minimal structure. 

\begin{figure}[htp]     \centering
    \includegraphics[width=14cm]{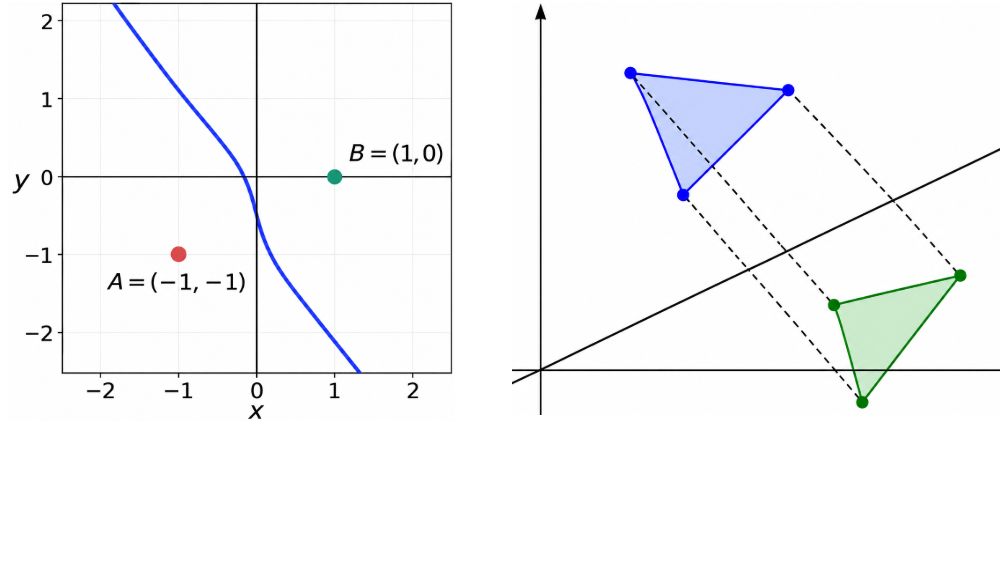}
\vspace{-22mm}    
    
    \caption{ Left: a non-line bisector under the $\ell_3$ norm. Right: an oblique reflection. }\label{fi:Oblique}
\end{figure}

\parag{Bisectors. } 
The bisector of two distinct points $p, q \in \R^2$ is $\B(p, q) = \{r \in \R^2 : \|p - r\| = \|q - r\|\}.$
In non-Euclidean norms, bisectors are not necessarily lines. 
For example, see the left side of Figure \ref{fi:Oblique}.
However, with strictly convex norms, every bisector is homeomorphic to a line (see \cite{Horvath00}).
We now characterize the line bisectors.

An \emph{oblique reflection} is a generalization of standard reflections. 
Such a reflection is still across a line, but the direction of reflection is not necessarily perpendicular to the line.
See the right side of Figure \ref{fi:Oblique}.
In particular, let $\ell_1,\ell_2$ be distinct lines through the origin.
Every vector $v\in \RR^2$ can be uniquely written as $v=u_1+u_2$, where $u_1\in\ell_1$ and $u_2\in \ell_2$.
The oblique reflection across $\ell_1$ in direction $\ell_2$ is $R(u_1+u_2) = u_1-u_2$.

With respect to the preceding paragraph, we note that (i) Every $r\in \ell_1$ satisfies $R(r)=r$; (ii) Every $r\in \ell_2$ satisfies $R(r)=-r$;
(iii) The composition $R\circ R$ is the identity; (iv) $R$ may not preserve Euclidean distances and angles.

An \emph{affine symmetry} of a curve $\gamma$ in $\RR^2$ is an invertible affine transformation that takes $\gamma$ to itself. 
For brevity, we will always write \emph{symmetry} when referring to an affine symmetry.

\begin{lemma} \label{whenline}
A bisector $\B(p, q)$ is a line if and only if the unit circle admits an oblique reflectional symmetry that preserves $p - q.$
\end{lemma}
\begin{proof}
We translate $\RR^2$ so that $\frac{p + q}{2} = 0$.
This leads to $q = -p$, which in turn implies that 
\[ \B(p, q) = \B(p, -p) = \{v \in \R^2\ :\ \|v - p\| = \|v + p\|\}. \]

First, we assume that $\B(p,-p)$ is a line. 
Let $v\in \RR^2$ be a vector in the direction of that line. 
By definition, for every $\lambda\in \RR$, we have that$\|\lambda v - p \| = \|\lambda v + p\|$. 
We note that $p$ and $v$ have different directions.
Thus, every point in $\R^2$ can be uniquely written as $\lambda v + \mu p$ with $\lambda, \mu \in \RR$. 
Then, the transformation $T(\lambda v + \mu p)= - \lambda v + \mu p$ is an oblique reflection that takes $p-q=2p$ to itself.
By the absolute homogeneity of norms, 
\[ \|T(\mu p+ \lambda v )\| = \|\mu p- \lambda v\|=|\mu|\cdot \|(\lambda/\mu)\cdot v-p\|=|\mu|\cdot \|(\lambda/\mu)\cdot v+p\|=\|\mu p+ \lambda v\|. \]
That is, the norm $\|\cdot\|$ is invariant under $T$.
In particular, $T$ is a symmetry of the unit circle.
 
For the other direction of the claim, we assume that the unit circle admits an oblique reflectional symmetry $T$ that preserves $p - q$. 
By the absolute homogeneity of a norm, $T$ is also a symmetry of any scaling of the unit circle.
In other words, $T$ preserves the length of any vector (formally, $T$ is an isometry).
Since $T$ is an oblique reflection about a line through the origin, there exists a nonzero $v \in \R^2$ that satisfies $T(v) = -v$.
Since $T$ is an isometry, every $\lambda \in \R$ satisfies 
\[ \|p + \lambda v\| = \|R(p + \lambda v)\| = \|p - \lambda v\|. \]

The above implies that $\lambda v \in \B(p, -p)$, for every $\lambda \in \RR$.
That is, $\B(p, -p)$ is the line through the origin of direction $v$.
\end{proof}

A planar convex body with infinitely many oblique reflection symmetries must be an ellipse (for example, see \cite[Proposition 2.6]{Dewar20}).
Thus, with non-ellipsoidal norms, line bisectors only have finally many directions.

\begin{corollary} \label{co:LineDirections}
For every non-ellipsoidal norm, there a constant $c$ such that all line bisectors have one of $c$ directions.
\end{corollary}

We will rely on another symmetry, which all bisectors have.

\begin{claim} \label{180_rot}
A bisector $\B(p,q)$ has a rotational symmetry of 180 degrees about the midpoint $\frac{p + q}{2}$.
\end{claim}
\begin{proof}
We translate the plane so that $\frac{p + q}{2} = (0,0)$.
We now have that $q = -p.$ 
It remains to show that every $r\in \B(p, q)$ satisfies $- r \in \B(p, q).$ 
Fix such $r \in \B(p, q)$.
By definition, $\|r - p\| = \|r + p\|$, which in turn implies that 
\[ \|(-r) - p\| = \|r + p\| = \|r - p\| = \|(-r) + p\|.\] 
We conclude that $-r \in \B(p, q)$, as required.
\end{proof}

\section{Bisectors in $o$-minimal Norms} \label{sec:bisectors}

Section \ref{sec:prelim} contained a few first properties of bisectors. 
In the current section, we study more involved properties of bisectors, some in the context of $o$-minimal norms. 

For the rest of this work, we consider only non-ellipsoidal strictly convex norms. 
Note that Theorem \ref{th:o-minimal} and Theorem \ref{th:semi-alg} consider only strictly convex norms. 
In addition, the Guth--Katz result \cite{GK15} and Claim \ref{cl:ellipsoidal} lead to a stronger bound for ellipsoidal norms. 
That is, Theorem \ref{th:o-minimal} and Theorem \ref{th:semi-alg} immediately hold for such norms.

\begin{figure}[htp]     \centering
    \includegraphics[width=6cm]{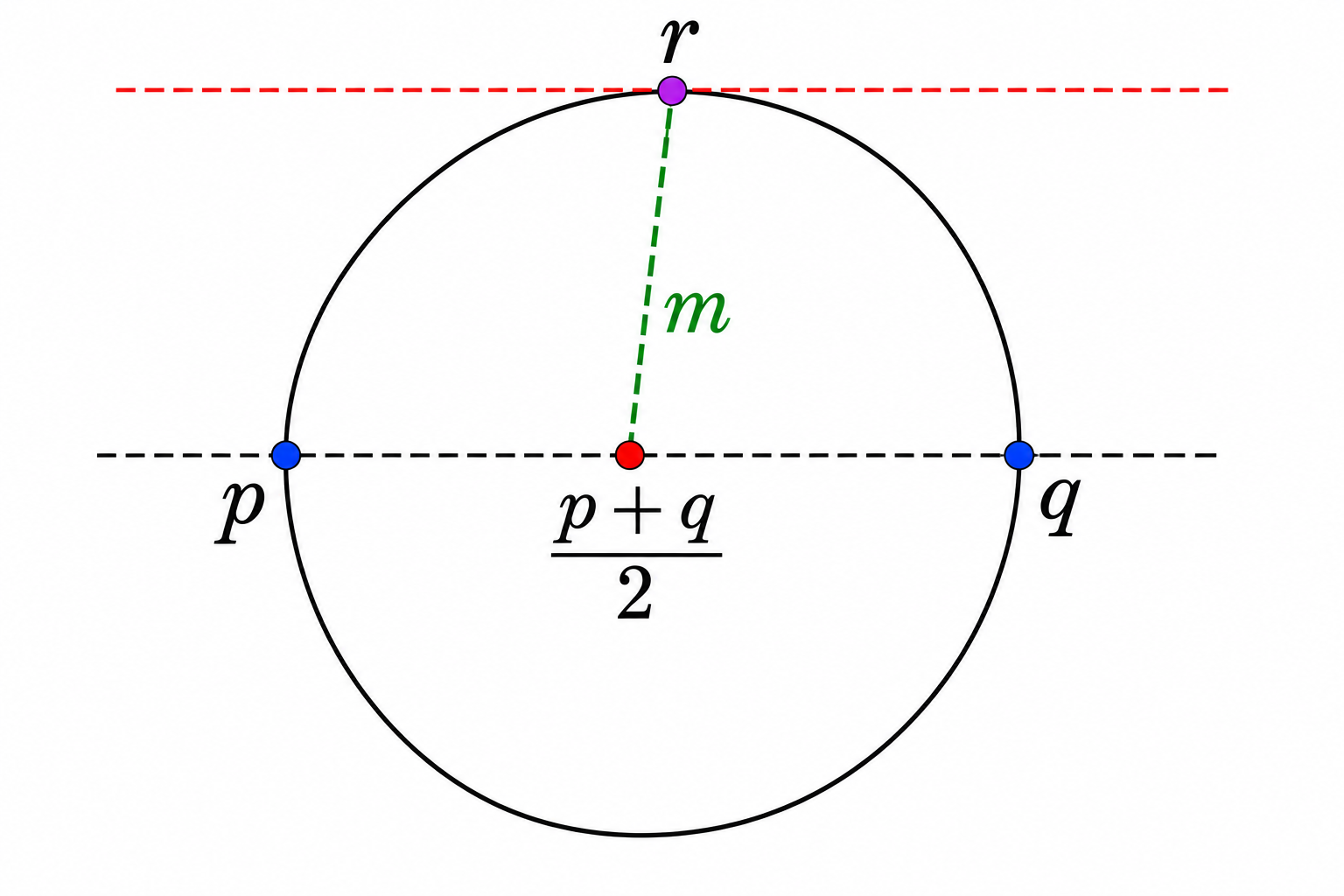}
\vspace{-5mm}    
    
    \caption{We consider the slope of the line between $r$ and the midpoint $(p+q)/2$.}\label{fi:TangentDirection}
\end{figure}

In the following, a \emph{circle} is the set of points of a specific distance from a center point. 
In particular, this is not a circle in the Euclidean sense.
The \emph{supporting line} of a point $r$ on a convex curve $\gamma$ is a line that contains $r$ and defines a closed half-plane that contains $\gamma$. 
See Figure \ref{fi:TangentDirection}.
For the following lemma, see Theorem 1 and Remark 1 of \cite{CMR96}.

\begin{lemma} \label{le:asymptoticLine}
For points $p,q\in \RR^2$, let $C$ be the circle centered at the midpoint $(p + q)/2$ and with radius $\|p - q\|/2$. 
Let $r$ be a point of $C$ with a supporting line parallel to the line through $p$ and $q$. 
Let $m$ be the slope of the line between the midpoint $(p + q)/2$ and $r$. (See Figure \ref{fi:TangentDirection}.)\\[2mm]
(a) If $B(p, q)$ has an asymptotic line, then the slope of this line is $m$. \\[2mm]
(b) If $B(p, q)$ does not have an asymptotic line, then $B(p, q)$ oscillates between two lines with slope $m$. That is, the bisector approaches each of the two lines infinitely many times.
\end{lemma}

Lemma \ref{le:asymptoticLine} implies the following.

\begin{lemma}\label{parallelcase}
If $p_1 - q_1$ and $p_2 - q_2$ have different directions, then $\B(p_1, q_1) \neq \B(p_2, q_2)$.
\end{lemma}
\begin{proof}
Let $C_1$ and $C_2$ be the circles defined in the statement of Lemma \ref{le:asymptoticLine}.
Let $m_1$ and $m_2$ be the slopes described in that lemma. 
These slopes are invariant under translating and uniformly scaling the plane. 
Thus, without loss of generality, we may assume that $(p_1+q_1)/2=(p_2+q_2)/2=(0,0)$ and that  $\|\frac{p_1 - q_1}{2}\| = \|\frac{p_2 - q_2}{2}\| = 1$.
In other words, $C_1=C_2$ is the unit circle around the origin.
We denote this circle as $C$.

Let $r_1$ and $r_2$ be the points described in Lemma \ref{le:asymptoticLine}.
Then $m_i$ is the slope of the line through the origin and $r_i$.
By the strict convexity and symmetry of the norm, we have that $m_1=m_2$ if and only if $r_1 = \pm r_2$.

Since the norm is strictly convex, no two points on $C$ have the same supporting line.
Thus, if $p_1 - q_1$ and $p_2 - q_2$ have different directions, then $r_1 \neq \pm r_2$.
This in turn implies that $m_1\neq m_2$.
By Lemma \ref{le:asymptoticLine}, we conclude that $\B(p_1, q_1) \neq \B(p_2, q_2)$.
\end{proof}

We now begin to study properties of bisectors under $o$-minimal norms.

\begin{lemma} \label{le:DefinableBisectors}
Consider a norm $\|\cdot\|$ that is definable in an o-minimal structure $\M$.
Then all bisectors under $\|\cdot\|$ are definable in $\M$.
\end{lemma}
\begin{proof}
Consider two points $p_1 = (x_1,y_1)$ and $p_2 = (x_2,y_2)$.
By the norm definition, when thinking of $x_1,y_1,x_2,y_2$ as variables, $(p_1, p_2) \mapsto \|p_1 - p_2\|$ is definable in $\M$. 
Then, for any fixed point $u\in \RR^2$, the restriction $p_1 \mapsto \|p_1 - u\|$ is definable as well. 
Indeed, this restriction is an intersection with a variety, followed by a coordinate-wise projection.
This in turn implies that the zero set of 
$x \mapsto \|x - u\| - \|x - v\|$ is definable.
This zero set is the bisector $\B(u,v)$.
\end{proof}

We are now ready to prove Theorem \ref{th:LineBisectors}.
We first recall the statement of this result.
\vspace{2mm}

\noindent {\bf Theorem \ref{th:LineBisectors}.}
\emph{Let $\|\cdot\|$ be a strictly convex $o$-minimal norm. 
If two pairs of points $(p_1,q_1) \neq (p_2,q_2)$ satisfy $\B(p_1, q_1) = \B(p_2, q_2)$, then this bisector is a line.}

\begin{proof}
Let $\gamma = \B(p_1, q_1) = \B(p_2, q_2).$ 
By Claim \ref{180_rot}, the bisector $\gamma$ has a rotational symmetry of 180 degrees about both midpoints $(p_1 + q_1)/2$ and $(p_2 + q_2)/2$.

We first consider the case where the two above midpoints are distinct.
In this case, the two rotational symmetries are distinct.
Combining the two distinct rotations by 180 degrees always yields a translation that is not the identity.
That is, $\gamma$ has a translation symmetry. 
Let $v$ be the vector of this translation. 
We arbitrarily pick a point $p\in \gamma$.
By the above, $\gamma$ contains the infinite set $\{p + nv : n \in \Z\}$.
By Lemma \ref{le:components}, $\gamma$ contains the line $\{p + av : a \in \RR\}$.
Since all bisectors are homeomorphic to lines, $\gamma$ is the above line.

It remains to consider the case where $(p_1 + q_1)/2=(p_2 + q_2)/2$.
Since translations do not affect the shape of a bisector, we may assume that the common midpoint is the origin. 
We now have that $q_1=-p_1$ and that $q_2=-p_2$.
By Lemma \ref{parallelcase}, $p_1-q_1$ and $p_2-q_2$ have the same direction.
That is, there exists $s\in \RR$ such that
\[ (p_1,-p_1) = (p_1,q_1) = s\cdot (p_2,q_2) = s\cdot (p_2,-p_2). \] 

We claim that $\gamma$ has a dilation symmetry.
Indeed, for any point $r\in \gamma$, we have that
\[ \|r - p_1\| = \|r + p_1\| \quad \Rightarrow \quad
    \|sr - sp_1\| = \|sr + sp_1\| \quad \Rightarrow \quad
    \|sr - p_2\| = \|sr + p_2\|. \]
That is, $sr \in \gamma$, as claimed.

We now fix a specific $r\in \gamma$ that is not the origin.
By the above, $\gamma$ contains the infinite set $\{s^n r : n \in \Z\}$.
As before, this implies that $\gamma$ contains the line $\{a r : a \in \RR\}$.
Since all bisectors are homeomorphic to lines, $\gamma$ is the above line.
\end{proof}

We also require the following property of definable bisectors.

\begin{lemma} \label{le:BisectorIntersection}
Let $\|\cdot\|$ be a strictly convex $o$-minimal norm, for which every two non-line bisectors intersect finitely many times. Then there exists $N$ such that every two non-line bisectors intersect at most $N$ times.
\end{lemma}
\begin{proof}
Let $\M$ be an $o$-minimal structure for which the norm is definable.
We consider the family of bisector intersections
\[ \left\{\B(p_1,q_1)\cap \B(p_2,q_2)\ :\ p_1,q_1,p_2,q_2\in \RR^2,\ p_1\neq q_1,\ p_2\neq q_2,\ (p_1,q_1)\neq (p_2,q_2) \right\}. \]
This is definable family of sets, as defined in Section \ref{sec:prelim}. 
It is an eight-dimensional family, depending on the eight coordinates of four points in $\RR^2$.

Consider two pairs of points $(p_1,q_1)\neq(p_2,q_2)$, such that neither $\B(p_1,q_1)$ nor $\B(p_2,q_2)$ is a line. 
By the assumption of the lemma, $\B(p_1,q_1)\cap \B(p_2,q_2)$ is finite.
Thus, Lemma \ref{le:components} implies that there exists an $N$ such that $\B(p_1,q_1)\cap \B(p_2,q_2)$ intersect in at most $N$ points.
\end{proof}

\section{Theorem \ref{th:o-minimal} and its applications} \label{sec:o-min}

In this section, we prove Theorem \ref{th:o-minimal} and its applications, Corollary \ref{co:lp} and Corollary \ref{co:o-minimal+analytic}.
These proofs build upon Garibaldi's work \cite{Garibaldi04} and its extension in \cite{Polymath20}.
We begin by collecting a few tools. 

\begin{lemma} \label{le:incidencebound}
    Let $\|\cdot\|$ be a strictly convex $o$-minimal norm for which every two non-line bisectors have a finite intersection. 
Let $\pts$ be a set of $m$ points in $\RR^2$ and let $\curves$ be a set of $n$ non-line bisectors. 
Then \[
    I(\pts, \curves) = O \left(m^{4/7} n^{6/7} + m + n \right).
    \]
\end{lemma}
\begin{proof}
By Lemma \ref{le:DefinableBisectors}, the elements of $\curves$ are definable under $\M$. 
We observe that a bisector $\B(p,q)$ is defined by the four coordinates of $p$ and $q$.
That is, the elements of $\curves$ belong to a four-dimensional family (as defined in Section \ref{sec:prelim}).
By Lemma \ref{le:BisectorIntersection}, there exists a constant $r$ such that two elements of $\curves$ intersect in at most $r-1$ points.
This implies that the incidence graph of $\pts\times\curves$ does not contain a $K_{r,r}$.
We may thus apply Lemma \ref{ominlemma} with $d=4$, which 
implies the assertion of the current lemma. 
\end{proof}

We recall that a \emph{multigraph} may have multiple edges between the same two vertices. 
The \emph{multiplicity} of a multigraph is the maximum number of edges that share the same two endpoints.
The \emph{crossing number} of a graph $G$ is the minimum number of pairs of intersecting edges in any planar drawing of $G$.
We denote this number as $cr(G)$.
For additional information about crossing numbers and the following lemma, see \cite{Szekely97}.
    
\begin{lemma} \label{CLfM}
Every multigraph $G = (V, E)$ with multiplicity $k$ and with $|E| > 5k \cdot |V|$ satisfies that  
\[   cr(G) = \Omega\left(\frac{|E|^3}{k \cdot |V|^2} \right).
        \]
    \end{lemma}

In the following, whenever we mention a circle, this is not a Euclidean circle. Instead, it is the set of points at a fixed distance from a given center, under the current norm. 
We are now ready to prove Theorem \ref{th:o-minimal}.
We first recall the statement of this result.
\vspace{2mm}

\noindent {\bf Theorem \ref{th:o-minimal}.}
\emph{Let $\|\cdot\|$ be a strictly convex $o$-minimal norm for which every two non-line bisectors have a finite intersection. Then, $d_{\|\cdot\|}(n) = \Omega(n^{5/6}).$ }

\begin{proof}
Let $\pts$ be a set of $n$ points and let $D$ be the set of distances spanned by $\pts$.
We set $t=|D|$.
Let $\curves$ be the set of all circles centered at a point of $\pts$ and with a radius from $D$.
That is, $\curves$ is a set of $nt$ circles.
For each $p\in \pts$, the $t$ circles centered at $p$ contain all points of $\pts\setminus\{p\}$. 
This implies that
\[  I(\pts, \curves) = n(n - 1) = \Theta(n^2). \]

We define a multigraph $G = (V, E)$, as follows.  
There is a vertex for each point of $\pts$.
For each circle $\gamma\in \curves$ that is incident to at least three points of $\pts$, we add an edge between every two vertices that appear consecutively on $\gamma$ (that is, with no third point between them).

We may assume that $t=O(n/\log n)$, since otherwise, the proof of the theorem is complete.
This implies that $|\curves|=O(n^2/\log n)$.
We remove from $\curves$ all circles that contain at most two points of $\pts$. 
By the above, after this step, we still have that $I(\pts,\curves)=\Theta(n^2)$.

Consider a remaining circle $\gamma\in\curves$ that is incident to $k\ge 3$ points of $\pts$.
By definition, $\gamma$ corresponds to $k$ edges of $G$. 
Thus, $|E|=I(\pts,\curves)=\Theta(n^2)$.

We can draw $G$ in $\RR^2$ by drawing the points of $\pts$ and circles of $\curves$ as given. 
Since each pair of circles intersects at most twice, the number of crossings in this drawing is at most $2\binom{|\curves|}{2}=O(n^2t^2)$.
This implies that $cr(G) = O(n^2 t^2)$.

We wish to apply Lemma \ref{CLfM}, to obtain a lower bound for $cr(G)$.
However, the multiplicity of $G$ may be too high for this to be helpful. 
To address this issue, we first discard from $G$ edges with high multiplicity. 

\parag{Removing high multiplicity edges.}
Consider an edge $(p,q)\in E$ that originates from $p$ and $q$ appearing consecutively on a circle around a point $r \in \pts$. 
By definition, $r\in \B(p, q)$. 
We remove from $G$ all such edges for which $\B(p, q)$ is a line.
By Corollary \ref{co:LineDirections}, all line bisectors have one of $O(1)$ specific slopes.
That is, every point $r\in \pts$ forms $O(1)$ incidences with line bisectors. 
For such a line, every element of $D$ corresponds to at most $n/2$ pairs $(p,q)$ that have this line as their bisector.
Indeed, every point of $\pts$ participates in at most one such pair.
We conclude that $O(nt)=O(n^2/\log n)$ edges originate from line bisectors.
After removing those, $G$ has $\Theta(n^2)$ remaining edges. 

Now, $G$ contains only edges $(p,q)$ where $\B(p,q)$ is a non-line bisectors. 
The multiplicity of such an edge is at most the number of points of $\pts$ on $\B(p,q)$.
Let $\curves_j$ be the set of non-line bisectors that contain at least $2^j$ points of $\pts$, but less than $2^{j+1}$.
By the definition of $\curves_j$ and Lemma \ref{le:incidencebound}, we have that
\[ 2^j \cdot |\curves_j|\le I(\pts,\curves_j) = O\left(n^{4/7}|\curves_j|^{6/7}+n+|\curves_j|\right). \] 

In the right-hand side of the above bound, when $j$ is larger than some fixed constant, the $+|\Gamma_j|$ term cannot dominate.
For every such $j$, rearranging the above leads to
\[ |\Gamma_j| = O\left(\frac{n^4}{2^{7j}}+\frac{n}{2^j}\right).\]

In the above bound, the first term dominates when $2^j = O(\sqrt{n})$.
For $k$ that will be set below, the number of edges in $G$ of multiplicity at least $k$ is smaller than
\begin{align*}
\sum_{j=\log k}^{\log n} |\curves_j|\cdot 2^{j+1} &= O\left(\sum_{j=\log k}^{0.5\cdot \log n} \frac{n^4}{2^{7j}}\cdot 2^j+\sum_{j=0.5\cdot \log n}^{\log n} \frac{n}{2^j}\cdot 2^j\right) \\[2mm]
&= O\left(\sum_{j=\log k}^{0.5\cdot \log n} \frac{n^4}{2^{6j}}+\sum_{j=0.5\cdot \log n}^{\log n} n\right) = O\left(\frac{n^4}{k^6}+n\log n\right).
\end{align*}

We set $k=cn^{1/3}$ for a constant $c$.
This implies that the number of edges of multiplicity at least $k$ is $O(n^2/c^6)$.
We recall that $G$ has $\Theta(n^2)$ edges. 
When $c$ is a sufficiently large constant, after removing all edges of multiplicity at least $k$, there are still $\Theta(n^2)$ edges in $G$.         
We perform this removal, obtaining that the multiplicity of $G$ is now smaller than $cn^{1/3}$.

\parag{Wrapping up.} 
With the above bound on the multiplicity of $G$, we may now apply Lemma \ref{CLfM}.
It is not difficult to verify that the lemma assumption $|E|>5 c n^{1/3} \cdot n$ indeed holds.
This leads to 
\[  cr(G) = \Omega\left(\frac{n^6}{cn^{1/3}\cdot n^2} \right) = \Omega(n^{11/3}). \]
    
We recall the above bound $cr(G)=O(n^2t^2)$.
Since removing edges cannot increase the number of crossings, it still holds. 
Combining these upper and lower bounds for $cr(G)$ implies that $n^{11/3} = O(n^2t^2)$.
Rearranging leads to $t=\Omega(n^{5/6})$, as required.
\end{proof}

Next, we move to show how Corollary \ref{co:lp} is obtained by applying Theorem \ref{th:o-minimal}.
We first recall the statement of this result.
\vspace{2mm}

\noindent {\bf Corollary \ref{co:lp}.}
\emph{   For all $p > 1,$ we have that $d_{\ell_p}(n) = \Omega(n^{5/6})$. }

\begin{proof}
The operations of taking real powers and absolute values are part of the $o$-minimal structure $\RR_{exp}$, as defined in Section \ref{sec:prelim} (for additional information, see for example \cite{Wilkie07}). 
That is, the $\ell_p$ norm $\|v\|_{\ell_p} = (|v_x|^p + |v_y|^p)^{1/p}$ is definable in $\RR_{exp}$.  
Garibaldi \cite{Garibaldi04} proved that, under $\ell_p$ norms, non-line bisectors have a finite intersection. 
We may thus apply Theorem \ref{th:o-minimal}.
\end{proof}

It remains to prove Corollary \ref{co:o-minimal+analytic}. We first recall the statement of this result.
\vspace{2mm}

\noindent {\bf Corollary \ref{co:o-minimal+analytic}.}
\emph{Let $\|\cdot\|$ be a strictly convex $o$-minimal norm with analytic bisectors. Then} $d_{\|\cdot\|}(n) = \Omega(n^{5/6})$.

\begin{proof}
By V\"ais\"al\"a \cite{Vaisala13}, for every non-line bisector $\B(p,q)$ there exists $f:\RR\to \RR$ such that  $\B(p,q)$ can be parameterized as $y=f(x)$.
By the assumption of the corollary, $f$ is analytic.

Consider two non-line bisectors that are parameterized by $f_1,f_2:\RR\to \RR$, as described above.
By Theorem \ref{th:LineBisectors}, the two bisectors are not identical.
The \emph{extension theorem} states that, if $f_1$ and $f_2$ have the same values on an open interval, then the two functions have the same values everywhere (for example, see \cite[Corollary 1.2.5]{KP02}).

By the above, $f_1$ and $f_2$ may not have the same values over an entire open interval.
Then, Lemma \ref{le:components} implies that the two bisectors have a finite intersection. 
We may thus apply Theorem \ref{th:o-minimal}.
\end{proof}

\section{Theorem \ref{th:semi-alg}} \label{sec:semi-alg}

In this section, we prove Theorem \ref{th:semi-alg}.
We begin with a semi-algebraic analog of Lemma \ref{le:incidencebound}.

\begin{lemma} \label{incidencebound}
Let $\|\cdot\|$ be a strictly convex semi-algebraic norm for which every two non-line bisectors have a finite intersection.
Let $\pts$ be a set of $m$ points in $\RR^2$ and let $\curves$ be a set of $n$ non-line bisectors.  
Then, for every $\eps>0$, we have that
\[ I(\pts, \curves) = O(m^{1/2} n^{7/8 + \epsilon} 
    + n^{2/3}m^{2/3} + m + n).
    \]
\end{lemma}
\begin{proof}
Each bisector is defined by two points, which consist of four coordinates. 
This implies that the elements of $\curves$ belong to a four-dimensional semi-algebraic family. 
By the assumption of the current lemma, two bisectors cannot share a one-dimensional component.

By Lemma \ref{le:BisectorIntersection}, there exists a constant $r$ such that two elements of $\curves$ intersect in at most $r-1$ points. 
We may thus apply Lemma \ref{semialglemma} with $s=4$.
This immediately leads to the statement of the current lemma. 
\end{proof}

We are now ready to prove Theorem \ref{th:semi-alg}. 
We first recall the statement of this result. 
\vspace{2mm}

\noindent {\bf Theorem \ref{th:semi-alg}.}
    \emph{Let $\|\cdot\|$ be a strictly convex semi-algebraic norm where bisectors intersect finitely many times. Then all $\eps > 0$ satisfy that $d_{\|\cdot\|}(n) = \Omega(n^{6/7 - \epsilon})$. }

\begin{proof}
We start the current proof in exactly the same way as the proof of Theorem \ref{th:o-minimal}.
Since semi-algebraic sets are also definable in an $o$-minimal structure, all steps off that proof are also valid here. 

For brevity, we do not repeat the identical beginning of the proof.
We define $t$, $D$, and the multigraph $G$ as before. 
We also remove from $G$ edges that correspond to line bisectors and obtain the bounds $cr(G)=(n^2t^2)$ and $|E|=\Theta(n^2)$, as before.

Let $\curves_j$ be the set of non-line bisectors from $\curves$ that are incident to at least $2^j$ points of $\pts$ and less than $2^{j+1}$ such points.  
We apply Lemma \ref{incidencebound} with each $\curves_j$. 
We do so with $\eps'>0$ whose value will be set below, to obtain that
\[ 2^j |\curves_j| \le I(\pts,\curves_j) = O\left(n^{1/2} |\curves_j|^{7/8 + \epsilon'} 
    + n^{2/3}|\curves_j|^{2/3} + n + |\curves_j|\right).\]

When $j$ is larger than a specific constant, term term $+|\curves_j|$ cannot dominate the above bound.
Rearranging the rest of the bound and taking $\eps'$ to be sufficiently small with respect to $\eps$ leads to 
\[ 
|\curves_j|= O\left(\frac{n^{4/(1-8\eps')}}{2^{8j/(1-8\eps')}}+\frac{n^2}{2^{3j}}+\frac{n}{2^j}\right) = O\left(\frac{n^{4+2\eps}}{2^{8j}}+\frac{n^2}{2^{3j}}+\frac{n}{2^j}\right).\]

In the above bound, the first term dominates when $j \le  \frac{2 + 2\eps}{5}\log n$, the second term dominates when $\frac{2 + 2\eps}{5}\log n < j \le \frac{1}{2} \log n,$ and the third when $j >\frac{1}{2} \log n$.
Thus, the number of remaining edges with multiplicity at least $k$ is smaller than
\begin{align*}
    \sum_{j=\log k}^{\log n} |\curves_j|\cdot 2^{j+1} &= O \left( 
    \sum_{j = \log k}^{\frac{2 + 2\eps}{5} \log n} \frac{n^{4 + 2\eps}}{2^{8j}} \cdot 2^j
    + \sum_{j = \frac{2 + 2\eps}{5} \log n}^{0.5\cdot \log n} \frac{n^{2}}{2^{3j}} \cdot 2^j
    + \sum_{0.5\cdot \log n}^{\log n} \frac{n}{2^j} \cdot 2^j
    \right)
   \\[2mm]
    &= O \left( 
     \frac{n^{4 + 2\eps}}{k^7} 
    + \frac{n^{2}}{(n^{(2+2\eps)/5})^2}     + n \log n \right) = O \left( 
     \frac{n^{4 + 2\eps}}{k^7} 
    + n^{(6-4\eps)/5}     \right).
\end{align*}
    
We set $k=n^{2/7+2\eps}$ and remove all edges with multiplicity at least $k$.
By the above, the number of removed edges is $O(n^{2-12\eps})$, so $G$ remains with $\Theta(n^2)$ edges.
Now, the multiplicity of $G$ is less than $n^{2/7+2\eps}$.
 
With the above multiplicity bound, we may now apply Lemma \ref{CLfM}.
It is not difficult to verify that the condition $|E|>5n^{2/7+2\eps}\cdot n$ indeed holds. 
This leads to 
\[ cr(G) = \Omega\left(\frac{n^6}{n^{2/7+2\eps}\cdot n^2}\right) = \Omega(n^{26/7-2\eps}).\]
 
We recall the above bound $cr(G)=O(n^2t^2)$. 
Since removing edges cannot increase the number of crossings, this bound still holds.
Combining the upper and lower bounds for $cr(G)$ implies that $n^{26/7-2\eps} = O(n^2t^2)$.
Rearranging this leads to $t=\Omega(n^{6/7-\eps})$, as required.
\end{proof}

To further advance the potato-based notation, we denote semi-algebraic norms where non-line bisectors have finite intersections as \emph{ube norms}. 
Ube norms are indeed a subset of sweet potato norms, since semi-algebraic sets are definable in every $o$-minimal structure.

\end{document}